\documentclass[11pt]{article}
\usepackage{latexsym}
\usepackage{amsfonts}
\usepackage{mathrsfs,amsthm,amsmath}
\usepackage{color}
\usepackage{todonotes}
\newtheorem{theorem}{Theorem}[section]
\newtheorem{proposition}[theorem]{Proposition}

\newtheorem{assumption}{Assumption}[section]
\newtheorem{remark}{Remark}[section]

\newtheorem{condition}{Condition}[section]
\begin{document}
\title{Noncentral moderate deviations for time-changed multivariate Lévy processes with linear combinations
of inverse stable subordinators\thanks{
C.M. acknowledges the support of MUR Excellence Department Project awarded to the Department of Mathematics, 
University of Rome Tor Vergata (CUP E83C23000330006), and by INdAM-GNAMPA.}}
\author{Neha Gupta\thanks{Address: Department of Mathematics and Statistics,
Indian Institute of Technology Kanpur, Kanpur 208016, India. e-mail: \texttt{nehagupta32022@gmail.com}}
\and Claudio Macci\thanks{Address: Dipartimento di Matematica, Università di Roma Tor Vergata,
Via della Ricerca Scientifica, I-00133 Rome, Italy. e-mail: \texttt{macci@mat.uniroma2.it}}}
\maketitle
\begin{abstract}
	The term noncentral moderate deviations is used in the literature to mean a class of large
	deviation principles that, in some sense, fills the gap between the convergence in probability to
	a constant (governed by a reference large deviation principle) and a weak convergence to a non-
	Gaussian (and non-degenerating) distribution. Some noncentral moderate deviation results in the 
	literature concern time-changed univariate Lévy processes, where the time-changes are given by 
	inverse stable subordinators. In this paper we present analogue results for multivariate Lévy 
	processes; in particular the random time-changes are suitable linear combinations of independent 
	inverse stable subordinators.
\end{abstract}

\section{Introduction}
The theory of large deviations gives an asymptotic computation of small probabilities on exponential scale (see 
\cite{DemboZeitouni} as a reference of this topic). The term \emph{moderate deviations} is used for a class of large 
deviation principles which fills the gap between a convergence to a constant $x_0$ (at least in probability), and a 
weak convergence to a non-constant centered Gaussian random variable. The convergence to $x_0$ is governed by a 
\emph{reference large deviation principle} with speed $v_t$, and a rate function that uniquely vanishes at $x_0$. 
Moderate deviations provide a class of large deviation principles which depend on the choice of some positive scalings 
$\{a_t:t>0\}$ such that $a_t\to 0$ and $v_ta_t\to\infty$ (as $t\to\infty$), each one with speed $1/a_t$ (thus the speed
is slower than $v_t$).

The term \emph{noncentral moderate deviations} has been recently introduced when one has the situation described above,
but the weak convergence is towards a non-constant and \emph{non-Gaussian} distributed random variable. A recent reference
with some examples and several references is \cite{GiulianoMacci}. Several examples in the literature concern families of
real valued random variables; a multivariate example is presented in \cite{LeonenkoMacciPacchiarotti}.

In this paper we consider light-tailed $h$-variate Lévy processes 
$$\{(S_1(t),\ldots,S_h(t)):t\geq 0\}$$
with independent components (see Assumption \ref{claim:on-Levy-processes}); obviously the univariate case 
$h=1$ is also allowed.
Our aim is to present noncentral moderate deviation results for two classes of independent time-changes of 
$\{(S_1(t),\ldots,S_h(t)):t\geq 0\}$ in terms of two classes of linear combinations of (independent) inverse stable 
subordinators; see Conditions \ref{cond:*} and \ref{cond:**} for a rigorous description.

The interest of time-changed processes is motivated by their potential applications in different areas
(e.g., finance, physics, and queueing theory) where random fluctuations can have an important impact on the system
dynamics. Indeed random time-changes can add a layer of complexity and realism to models. Several examples of random
time-changes are given by nondecreasing Lévy processes (called \emph{subordinators}) or their inverses. Here we recall some
references in which some inverse of subordinators (or their mixtures) are studied:
\cite{CD} for inverse Gamma subordinators,
\cite{MS} for inverse stable subordinators,
\cite{LMSS} for inverse of mixtures of stable subordinators,
\cite{KV} for inverse tempered stable subordinators, and
\cite{GKL} for inverse of mixtures of tempered stable subordinators.
We also recall that the inverse stable subordinators have important connections with time-fractional processes.

We start by introducing the light-tailed $h$-variate Lévy processes $\{(S_1(t),\ldots,S_h(t)):t\geq 0\}$ with independent 
components; moreover we introduce some related notation.

\begin{assumption}\label{claim:on-Levy-processes}
	Let $\{S_1(t):t\geq 0\},\ldots,\{S_h(t):t\geq 0\}$ be independent real-valued Lévy processes;
	moreover they are assumed to be light-tailed, i.e. the functions
	$\kappa_{S_1},\ldots,\kappa_{S_h}$ defined by
	$$\kappa_{S_i}(\theta):=\log\mathbb{E}[e^{\theta S_i(1)}]\quad (\mbox{for all}\ i\in\{1,\ldots,h\}\ 
	\mbox{and}\ \theta\in\mathbb{R})$$
	are finite in a neighborhood of $\theta=0$, and therefore the vector
	$$m=(m_1,\ldots,m_h):=(\kappa_{S_1}^\prime(0),\ldots,\kappa_{S_h}^\prime(0))=(\mathbb{E}[S_1(1)],\ldots,\mathbb{E}[S_h(1)])$$
	is well-defined.
\end{assumption}

We have the following remarks. 

\begin{remark}\label{rem:degenerating-case-with-m=0}
	If $\kappa_{S_i}^{\prime\prime}(0)=0$ for some $i\in\{1,\ldots,h\}$, then $\{S_i(t):t\geq 0\}$ is a deterministic Lévy process
	and, with the notation in Assumption \ref{claim:on-Levy-processes}, we have $S_i(t)=m_it$. Thus we can say that, if 
	$m=(m_1,\ldots,m_h)=0\in\mathbb{R}^h$ and $\kappa_{S_i}^{\prime\prime}(0)=0$ for every $i\in\{1,\ldots,h\}$, then the processes
	in Conditions \ref{cond:*} and \ref{cond:**} below are identically equal to zero and all the results below will be trivial;
	moreover, in such a case, the weak convergence results in Propositions \ref{prop:weak-convergence*} and
	\ref{prop:weak-convergence**} are towards the constant random variable equal to $0\in\mathbb{R}^h$.
\end{remark}

\begin{remark}\label{rem:discussion}
	Here we present a brief discussion on Assumption \ref{claim:on-Levy-processes}.
	The independence of the processes $\{S_1(t):t\geq 0\},\ldots,\{S_h(t):t\geq 0\}$ allows to make calculations on each component 
	of the random vectors separately; therefore we can refer to the asymptotic behavior of the Mittag-Leffler functions involved, and
	we can refer to what happens for the univariate case. In our opinion this difficulty could be overcome at least under 
	suitable hypotheses, and this could be studied in a successive work.
	Moreover $\{S_1(t):t\geq 0\},\ldots,\{S_h(t):t\geq 0\}$ are assumed to be light-tailed, and this allows to apply the G\"artner 
	Ellis Theorem (Theorem \ref{th:GE} below). In general this problem would not be easy to face. For instance, as discussed in Section 
	4 in \cite{BeghinMacciSPL2022} (where the Lévy processes are compound Poisson processes), one could hope to find a connection 
	with the results on empirical means of i.i.d. semi-exponential distributed random variables (see the references cited therein).
\end{remark}

Now we present two classes of independent time-changes for $\{(S_1(t),\ldots,S_h(t)):t\geq 0\}$ studied in this paper.
In particular we refer to the concept of inverse stable subordinator recalled below (see just after 
eq. \eqref{eq:MGF-inverse-stable-sub}).
We remark that, actually, some of the positive coefficients $c_i$ in Condition \ref{cond:*} and $c_{ij}$ in Condition 
\ref{cond:**} could be equal to zero, and the statements of the results can be suitably modified (the details will 
be omitted).

\begin{condition}\label{cond:*}
	Let $c_0,c_1,\ldots,c_h>0$ be arbitrarily fixed. Then let
	$$\{(S_1(c_1L_1^{\nu_1}(t)+c_0L_0^{\nu_0}(t)),\ldots,S_h(c_hL_h^{\nu_h}(t)+c_0L_0^{\nu_0}(t))):t\geq 0\}$$
	be the stochastic process defined in terms of the following families of independent random variables:
\begin{itemize}
	\item $\{S_1(t):t\geq 0\},\ldots,\{S_h(t):t\geq 0\}$ are as in Assumption \ref{claim:on-Levy-processes};
	\item $\{L_0^{\nu_0}(t):t\geq 0\},\{L_1^{\nu_1}(t):t\geq 0\},\ldots,\{L_h^{\nu_h}(t):t\geq 0\}$ are independent
	inverse stable subordinators of parameters $\nu_0,\nu_1,\ldots,\nu_h\in(0,1)$, respectively.
\end{itemize}
\end{condition}

\begin{condition}\label{cond:**}
	Let $\{c_{ij}:i\in\{1,\ldots,h\},j\in\{1,\ldots,k\}\}$ be arbitrarily fixed positive 
	numbers. Then let
	$$\left\{\left(S_1\left(\sum_{j=1}^kc_{1j}L_j^\nu(t)\right),\ldots,
	S_h\left(\sum_{j=1}^kc_{hj}L_j^\nu(t)\right)\right):t\geq 0\right\}$$
	be the stochastic process defined in terms of the following families of independent random variables:
	\begin{itemize}
		\item $\{S_1(t):t\geq 0\},\ldots,\{S_h(t):t\geq 0\}$ are as in Assumption \ref{claim:on-Levy-processes};
		\item $\{L_1^\nu(t):t\geq 0\},\ldots,\{L_k^\nu(t):t\geq 0\}$ are independent
		inverse stable subordinators of parameter $\nu\in(0,1)$.
	\end{itemize}
\end{condition}

The results in this paper have some analogies with the ones for univariate time-changed Lévy processes in the 
literature: see e.g. \cite{IulianoMacciMeoli}, with some results for Lévy processes with bounded variation 
(and therefore for differences of two independent subordinators); see also \cite{BeghinMacciSPL2022} which concerns
compound Poisson processes only, and \cite{LeeMacci} which concerns Skellam processes (which are differences between 
independent Poisson processes). This paper is the first contribution with results for the multivariate case.

Under both Conditions \ref{cond:*} and \ref{cond:**} we present the following results.
\begin{itemize}
	\item A reference LDP with speed $v_t=t$ (which does not depend on $m$ in Assumption \ref{claim:on-Levy-processes}),
	and we have a convergence in probability to $x_0=0$, where $0\in\mathbb{R}^h$ is the null vector; see Proposition
	\ref{prop:LD*} under Condition \ref{cond:*}, and Proposition \ref{prop:LD**} under Condition \ref{cond:**}.
	\item Two weak convergence results (for $m=0$ and $m\neq 0$, where $0\in\mathbb{R}^h$ is the null vector) towards 
	non-constant and non-Gaussian distributed random variables; see Proposition \ref{prop:weak-convergence*} under
	Condition \ref{cond:*}, and Proposition \ref{prop:weak-convergence**} under Condition \ref{cond:**}.
	\item Two noncentral moderate deviation results (for $m=0$ and $m\neq 0$, again) concerning classes of LDPs which 
	depend on the choices of some positive scalings $\{a_t:t\geq 0\}$ such that
	\begin{equation}\label{eq:MDconditions}
		a_t\to 0\ \mbox{and}\ ta_t\to\infty\quad (\mbox{as}\ t\to\infty);
	\end{equation}
    see Proposition \ref{prop:ncMD*} under Condition \ref{cond:*}, and Proposition \ref{prop:ncMD**} under Condition 
    \ref{cond:**}.
\end{itemize}

We can say that the results in this paper have some features of the results for the univariate processes in the literature;
see \cite{IulianoMacciMeoli} (Section 3), or in \cite{BeghinMacciSPL2022} (for compound Poisson processes only). We mean that 
the reference LDP that does not depend on $m$ in Assumption \ref{claim:on-Levy-processes} while, for the other 
results, we have to distinguish the cases $m=0$ and $m\neq 0$. To this aim it is useful to introduce the notation
\begin{equation}\label{eq:exponents}
	\alpha_m(\nu):=\left\{\begin{array}{ll}
		1-\nu/2&\ \mbox{if}\ m=0\\
		1-\nu&\ \mbox{if}\ m\neq 0,
	\end{array}\right.\quad\mbox{for}\ \nu\in(0,1).
\end{equation}
Moreover, for the results concerning Condition \ref{cond:*}, we use the notation $x\vee y$ to denote 
the maximum between $x,y\in\mathbb{R}$.

We conclude with the outline of the paper. Section \ref{sec:preliminaries} is devoted to recall some preliminaries. 
In Section \ref{sec:results*} we prove the results under Condition \ref{cond:*}, and in Section 
\ref{sec:explicit-expressions} we present some cases in which it is possible to get an explicit expression for the rate 
function in Proposition \ref{prop:ncMD*} (concerning moderate deviations). In Section \ref{sec:results**} we prove the 
results under Condition \ref{cond:**}.

\section{Preliminaries}\label{sec:preliminaries}
In this section we recall some preliminaries on large deviations, on the Lévy processes in Assumption 
\ref{claim:on-Levy-processes}, and on the inverse the stable subordinator, together with the Mittag-Leffler function.

\subsection{Preliminaries on large deviations}
We start with some basic definitions (see e.g. \cite{DemboZeitouni}). In view of what follows we present definitions 
and results for families of $\mathbb{R}^b$-valued random variables $\{Z(t):t>0\}$, for some 
$b$, defined on the same probability space $(\Omega,\mathcal{F},P)$, where $t$ goes to infinity. Throughout this
paper we always have $b=h$.

A real-valued function $\{v_t:t>0\}$ such that $v_t\to\infty$ (as $t\to\infty$) is called a \emph{speed function}, 
and a lower semicontinuous function 
$I:\mathbb{R}^b\to[0,\infty]$ is called a \emph{rate function}. Then $\{Z(t):t>0\}$ satisfies the LDP with speed $v_t$ 
and a rate function $I$ if
$$\limsup_{t\to\infty}\frac{1}{v_t}\log P(Z(t)\in C)\leq-\inf_{x\in C}I(x)\quad\mbox{for all closed sets}\ C,$$
and
$$\liminf_{t\to\infty}\frac{1}{v_t}\log P(Z(t)\in O)\geq-\inf_{x\in O}I(x)\quad\mbox{for all open sets}\ O.$$
The rate function $I$ is said to be \emph{good} if, for every $\beta\geq 0$, the level set 
$\{x\in\mathbb{R}^b:I(x)\leq\beta\}$ is compact. We also recall the following known result (see e.g. Theorem
2.3.6(c) in \cite{DemboZeitouni}).

\begin{theorem}[G\"artner Ellis Theorem]\label{th:GE}
Assume that, for all $\theta\in\mathbb{R}^b$, there exists
$$\Lambda(\theta):=\lim_{t\to\infty}\frac{1}{v_t}\log\mathbb{E}\left[e^{v_t\langle\theta,Z(t)\rangle}\right]$$
as an extended real number (here $\langle\cdot,\cdot\rangle$ is the inner product in $\mathbb{R}^b$). Assume that 
$\theta=0\in\mathbb{R}^b$ belongs to the interior of the set 
$\mathcal{D}(\Lambda):=\{\theta\in\mathbb{R}^b:\Lambda(\theta)<\infty\}$; moreover assume that the function 
$\Lambda$ is essentially smooth and lower semi-continuous. Then the family of random variables $\{Z(t):t>0\}$
satisfies the LDP with good rate function $\Lambda^*$ defined by
$$\Lambda^*(x):=\sup_{\theta\in\mathbb{R}^b}\{\langle\theta,x\rangle-\Lambda(\theta)\}$$
(i.e. $\Lambda^*$ is the Legendre-Fenchel transform of $\Lambda$).
\end{theorem}

We also recall (see e.g. Definition 2.3.5 in \cite{DemboZeitouni}) that $\Lambda$ is essentially smooth
if the interior of $\mathcal{D}(\Lambda)$ is non-empty, the function $\Lambda$ is differentiable 
throughout the interior of $\mathcal{D}(\Lambda)$, and $\Lambda$ is steep, i.e. $\|\Lambda^\prime(\theta_n)\|\to\infty$
whenever $\theta_n$ is a sequence of points in the interior of $\mathcal{D}(\Lambda)$ which converge to 
a boundary point of $\mathcal{D}(\Lambda)$.

\subsection{Preliminaries on Lévy processes in Assumption \ref{claim:on-Levy-processes}}
Some well-known properties on Lévy processes will be used for the processes in Assumption \ref{claim:on-Levy-processes}.
Firstly, for every $t\geq 0$, we have
$$\mathbb{E}[e^{\theta S_i(t)}]=e^{t\kappa_{S_i}(\theta)}\quad (\mbox{for all}\ i\in\{1,\ldots,h\}\ 
\mbox{and}\ \theta\in\mathbb{R}).$$
Moreover, for every $i\in\{1,\ldots,h\}$, we consider the following Taylor formulas for the function
$\kappa_{S_i}$:
$$\kappa_{S_i}(\theta)=m_i\theta+o(\theta)\quad\mbox{as}\ \theta\to 0$$
and, if $m=0\in\mathbb{R}^h$, i.e. $m_1=\cdots=m_h=0$,
$$\kappa_{S_i}(\theta)=\kappa_{S_i}^{\prime\prime}(0)\frac{\theta^2}{2}+o(\theta^2)\quad\mbox{as}\ \theta\to 0.$$
We recall the discussion in Remark \ref{rem:degenerating-case-with-m=0} on the case $\kappa_{S_i}^{\prime\prime}(0)=0$.

\subsection{Preliminaries on inverse stable subordinators}
We start with the definition of the Mittag-Leffler function (see e.g. \cite{GorenfloKilbasMainardiRogosin}, eq. 
(3.1.1))
$$E_\alpha(x):=\sum_{k=0}^\infty\frac{x^k}{\Gamma(\alpha k+1)}.$$
It is known (see Proposition 3.6 in \cite{GorenfloKilbasMainardiRogosin} for the case $\alpha\in(0,2)$; indeed here
we consider the case $\alpha\in(0,1)$ only) that we have
\begin{equation}\label{eq:ML-asymptotics}
\left\{\begin{array}{l}
E_\alpha(x)\sim\frac{e^{x^{1/\alpha}}}{\alpha}\ \mbox{as}\ x\to\infty\\
\mbox{if}\ y<0,\ \mbox{then}\ \frac{1}{x}\log E_\alpha(xy)\to 0\ \mbox{as}\ x\to\infty.
\end{array}\right.
\end{equation}
Then, for the stable subordinator of order $\nu\in(0,1)$, denoted by $\{S_\nu(t):t\geq 0\}$, for all $t\geq 0$ we have
$$\mathbb{E}[e^{\theta S_\nu(t)}]=\left\{\begin{array}{ll}
	e^{-t(-\theta)^\nu}&\quad\mbox{if}\ \theta\leq 0\\
	\infty&\quad\mbox{if}\ \theta>0.
\end{array}\right.$$ 
Moreover, for the inverse stable subordinator $\{L_\nu(t):t\geq 0\}$ (defined by
$$L_\nu(t):=\inf\{s\geq 0:S_\nu(s)>t\},$$
where $\{S_\nu(t):t\geq 0\}$ is a stable subordinator of order $\nu\in(0,1)$), for all $t\geq 0$ we have 
\begin{equation}\label{eq:MGF-inverse-stable-sub}
\mathbb{E}[e^{\theta L_\nu(t)}]=E_\nu(\theta t^\nu)\quad \mbox{for all}\ \theta\in\mathbb{R}.
\end{equation}

\section{Results under Condition \ref{cond:*}}\label{sec:results*}
We start with the reference LDP and the result does not depend on $m$ in Assumption \ref{claim:on-Levy-processes}.
In particular one can check that $I_{\mathrm{LD}}(x_1,\ldots,x_h)=0$ if and only if $(x_1,\ldots,x_h)=(0,\ldots,0)$.

\begin{proposition}\label{prop:LD*}
	Assume that Condition \ref{cond:*} holds.
	Let $\Psi$ be the function defined by
	\begin{multline*}
		\Psi(\theta_1,\ldots,\theta_h):=\sum_{i=1}^h(c_i\kappa_{S_i}(\theta_i))^{1/\nu_i}1_{\{\kappa_{S_i}(\theta_i)\geq 0\}}\\
		+\left(c_0\sum_{i=1}^h\kappa_{S_i}(\theta_i)\right)^{1/\nu_0}1_{\{\sum_{i=1}^h\kappa_{S_i}(\theta_i)\geq 0\}}
		\quad(\mbox{for all}\ \theta_1,\ldots,\theta_h\geq 0),
	\end{multline*}
	and assume that it is an essentially smooth and lower semicontinuous function. Then
	$$\left\{\left(\frac{S_1(c_1L_1^{\nu_1}(t)+c_0L_0^{\nu_0}(t))}{t},\ldots,
	\frac{S_h(c_hL_h^{\nu_h}(t)+c_0L_0^{\nu_0}(t))}{t}\right):t>0\right\}$$
	satisfies the LDP with speed function $t$ and rate function $I_{\mathrm{LD}}$ defined by
	$$I_{\mathrm{LD}}(x_1,\ldots,x_h):=\sup_{(\theta_1,\ldots,\theta_h)\in\mathbb{R}^h}\left\{\sum_{i=1}^h\theta_ix_i
	-\Psi(\theta_1,\ldots,\theta_h)\right\}\quad(\mbox{for all}\ (x_1,\ldots,x_h)\in\mathbb{R}^h.$$
\end{proposition}
\begin{proof}
	The desired LDP can be derived by applying the G\"artner Ellis Theorem (i.e. Theorem \ref{th:GE}). In fact, for all
	$t>0$ (and for all $(\theta_1,\ldots,\theta_h)\in\mathbb{R}^h$), we have
	\begin{multline*}
		\mathbb{E}\left[e^{\sum_{i=1}^h\theta_iS_i(c_iL_i^{\nu_i}(t)+c_0L_0^{\nu_0}(t))}\right]
		=\mathbb{E}\left[\mathbb{E}\left[e^{\sum_{i=1}^h\theta_iS_i(r_i+r_0)}\right]_{r_0=c_0L_0^{\nu_0}(t),
		r_1=c_1L_1^{\nu_1}(t),\ldots,r_h=c_hL_h^{\nu_h}(t)}\right]\\
	    =\mathbb{E}\left[\prod_{i=1}^h\mathbb{E}\left[e^{\theta_iS_i(1)}\right]^{c_iL_i^{\nu_i}(t)+c_0L_0^{\nu_0}(t)}\right]
	    =\mathbb{E}\left[e^{\sum_{i=1}^h(c_iL_i^{\nu_i}(t)+c_0L_0^{\nu_0}(t))\kappa_{S_i}(\theta_i)}\right]\\
	    =\mathbb{E}\left[e^{\sum_{i=1}^hc_iL_i^{\nu_i}(t)\kappa_{S_i}(\theta_i)+
	    c_0L_0^{\nu_0}(t)\sum_{i=1}^h\kappa_{S_i}(\theta_i)}\right]
	\end{multline*}
	whence we obtain (by the independence of the involved inverse stable subordinators, and by the expressions of the
	involved moment generating functions given by \eqref{eq:MGF-inverse-stable-sub} with one among $\nu_0,\nu_1,\ldots,\nu_h$
	in place of $\nu$)
	\begin{multline*}
		\frac{1}{t}\log\mathbb{E}\left[e^{\sum_{i=1}^h\theta_iS_i(c_iL_i^{\nu_i}(t)+c_0L_0^{\nu_0}(t))}\right]
		=\sum_{i=1}^h\frac{1}{t}\log\mathbb{E}\left[e^{c_iL_i^{\nu_i}(t)\kappa_{S_i}(\theta_i)}\right]
		+\frac{1}{t}\log\mathbb{E}\left[e^{c_0L_0^{\nu_0}(t)\sum_{i=1}^h\kappa_{S_i}(\theta_i)}\right]\\
		=\sum_{i=1}^h\frac{1}{t}\log E_{\nu_i}(c_i\kappa_{S_i}(\theta_i)t^{\nu_i})
		+\frac{1}{t}\log E_{\nu_0}\left(c_0\sum_{i=1}^h\kappa_{S_i}(\theta_i)t^{\nu_0}\right).
	\end{multline*}
	Finally, if we take the limit as $t$ tends to infinity, we get
	$$\lim_{t\to\infty}\frac{1}{t}\log\mathbb{E}\left[e^{\sum_{i=1}^h\theta_iS_i(c_iL_i^{\nu_i}(t)+c_0L_0^{\nu_0}(t))}\right]
	=\Psi(\theta_1,\ldots,\theta_h)$$
	by taking into account eq. \eqref{eq:ML-asymptotics} with $\alpha=\nu_i$ for all $i\in\{1,\ldots,h\}$, and for $\alpha=\nu_0$.
\end{proof}

In the following results we have to distinguish the cases $m=0$ and $m\neq 0$, where $0\in\mathbb{R}^h$ is the null vector.
We start with the weak convergence result; in our setting the weak convergence can be proved by taking the limit of moment 
generating functions (instead of characteristic functions).

\begin{proposition}\label{prop:weak-convergence*}
	Assume that Condition \ref{cond:*} holds.	
	Let $\alpha_m(\nu)$ be defined in eq. \eqref{eq:exponents}. We have the following statements.
	\begin{itemize}
		\item If $m=0$, then $\left\{\left(t^{\alpha_m(\nu_0\vee\nu_i)}
		\frac{S_i(c_iL_i^{\nu_i}(t)+c_0L_0^{\nu_0}(t))}{t}\right)_{i\in\{1,\ldots,h\}}:t>0\right\}$ 
		converges weakly to
		$$\left(\sqrt{c_iL_i^{\nu_i}(1)\kappa_{S_i}^{\prime\prime}(0)1_{\{\nu_0\leq\nu_i\}}}Z_i+
		\sqrt{c_0L_0^{\nu_0}(1)\kappa_{S_i}^{\prime\prime}(0)1_{\{\nu_i\leq\nu_0\}}}\widehat{Z}_i\right)_{i\in\{1,\ldots,h\}},$$
		where $Z_1,\ldots,Z_h,\widehat{Z}_1,\ldots,\widehat{Z}_h$ are independent standard Normal distributed random 
		variables, and independent of all the other random variables.
		\item If $m\neq 0$, then $\left\{\left(t^{\alpha_m(\nu_0\vee\nu_i)}
		\frac{S_i(c_iL_i^{\nu_i}(t)+c_0L_0^{\nu_0}(t))}{t}\right)_{i\in\{1,\ldots,h\}}:t>0\right\}$ 
		converges weakly to
		$$\left(c_im_i1_{\{\nu_0\leq\nu_i\}}L_i^{\nu_i}(1)+c_0m_i1_{\{\nu_i\leq\nu_0\}}L_0^{\nu_0}(1)\right)_{i\in\{1,\ldots,h\}}.$$
	\end{itemize}
\end{proposition}
\begin{proof}
	In both cases $m=0$ and $m\neq 0$ we study suitable limits (as $t\to\infty$) of suitable moment generating functions,
	which can be expressed in terms of Mittag-Leffler functions (see \eqref{eq:MGF-inverse-stable-sub}, with one among 
	$\nu_0,\nu_1,\ldots,\nu_h$ in place of $\nu$). Moreover
	\begin{multline*}
	    \mathbb{E}\left[e^{\sum_{i=1}^h\theta_it^{\alpha_m(\nu_0\vee\nu_i)}\frac{S_i(c_iL_i^{\nu_i}(t)+c_0L_0^{\nu_0}(t))}{t}}\right]
	    =\mathbb{E}\left[e^{\sum_{i=1}^h\theta_i\frac{S_i(c_iL_i^{\nu_i}(t)+c_0L_0^{\nu_0}(t))}{t^{1-\alpha_m(\nu_0\vee\nu_i)}}}\right]\\
	    =\mathbb{E}\left[\mathbb{E}\left[e^{\sum_{i=1}^h\theta_i\frac{S_i(r_i+r_0)}{t^{\alpha_m(\nu_0\vee\nu_i)}}}
	    \right]_{r_0=c_0L_0^{\nu_0}(t),r_1=c_1L_1^{\nu_1}(t),\ldots,r_h=c_hL_h^{\nu_h}(t)}\right]\\
	    =\mathbb{E}\left[e^{\sum_{i=1}^h(c_iL_i^{\nu_i}(t)+c_0L_0^{\nu_0}(t))\kappa_{S_i}(\theta_i/t^{1-\alpha_m(\nu_0\vee\nu_i)})}\right]\\
	    =\prod_{i=1}^h\mathbb{E}\left[e^{c_i\kappa_{S_i}(\theta_i/t^{1-\alpha_m(\nu_0\vee\nu_i)})L_i^{\nu_i}(t)}\right]
	    \mathbb{E}\left[e^{c_0\sum_{i=1}^h\kappa_{S_i}(\theta_i/t^{1-\alpha_m(\nu_0\vee\nu_i)})L_0^{\nu_0}(t)}\right]\\ 
	    =\prod_{i=1}^hE_{\nu_i}\left(c_i\kappa_{S_i}(\theta_i/t^{1-\alpha_m(\nu_0\vee\nu_i)})t^{\nu_i}\right)
	    E_{\nu_0}\left(c_0\sum_{i=1}^h\kappa_{S_i}(\theta_i/t^{1-\alpha_m(\nu_0\vee\nu_i)})t^{\nu_0}\right).
    \end{multline*}
	
	We start with the case $m=0$. We have
	\begin{multline*}
		\mathbb{E}\left[e^{\sum_{i=1}^h\theta_it^{\alpha_m(\nu_0\vee\nu_i)}\frac{S_i(c_iL_i^{\nu_i}(t)+c_0L_0^{\nu_0}(t))}{t}}\right]\\
		=\prod_{i=1}^hE_{\nu_i}\left(c_i\left(\frac{\kappa_{S_i}^{\prime\prime}(0)}{2}\frac{\theta_i^2}{t^{\nu_0\vee\nu_i}}
		+o\left(\frac{1}{t^{\nu_0\vee\nu_i}}\right)\right)t^{\nu_i}\right)
		E_{\nu_0}\left(c_0\sum_{i=1}^h\left(\frac{\kappa_{S_i}^{\prime\prime}(0)}{2}\frac{\theta_i^2}{t^{\nu_0\vee\nu_i}}
		+o\left(\frac{1}{t^{\nu_0\vee\nu_i}}\right)\right)t^{\nu_0}\right)
	\end{multline*}
    and, by \eqref{eq:ML-asymptotics}, one can check that the limit of the moment generating function (as $t\to\infty$)
    is equal to
	\begin{multline*}
		\prod_{i=1}^hE_{\nu_i}\left(\frac{c_i\kappa_{S_i}^{\prime\prime}(0)}{2}1_{\{\nu_0\leq\nu_i\}}\theta_i^2\right)
		E_{\nu_0}\left(c_0\sum_{i=1}^h\frac{\kappa_{S_i}^{\prime\prime}(0)}{2}1_{\{\nu_i\leq\nu_0\}}\theta_i^2\right)\\
		=\prod_{i=1}^h\mathbb{E}\left[e^{\theta_i\sqrt{c_iL_i^{\nu_i}(1)\kappa_{S_i}^{\prime\prime}(0)1_{\{\nu_0\leq\nu_i\}}}Z_i}\right]
		\mathbb{E}\left[e^{\sum_{i=1}^h\theta_i\sqrt{c_0L_0^{\nu_0}(1)\kappa_{S_i}^{\prime\prime}(0)1_{\{\nu_i\leq\nu_0\}}}\widehat{Z}_i}\right]\\
		=\mathbb{E}\left[e^{\sum_{i=1}^h\theta_i\left(\sqrt{c_iL_i^{\nu_i}(1)\kappa_{S_i}^{\prime\prime}(0)1_{\{\nu_0\leq\nu_i\}}}Z_i+
		\sqrt{c_0L_0^{\nu_0}(1)\kappa_{S_i}^{\prime\prime}(0)1_{\{\nu_i\leq\nu_0\}}}\widehat{Z}_i\right)}\right].
	\end{multline*}
	
	Finally the case $m\neq 0$. We have
	\begin{multline*}
		\mathbb{E}\left[e^{\sum_{i=1}^h\theta_it^{\alpha_m(\nu_0\vee\nu_i)}\frac{S_i(c_iL_i^{\nu_i}(t)+c_0L_0^{\nu_0}(t))}{t}}\right]\\
		=\prod_{i=1}^hE_{\nu_i}\left(c_i\left(\kappa_{S_i}^\prime(0)\frac{\theta_i}{t^{\nu_0\vee\nu_i}}
		+o\left(\frac{1}{t^{\nu_0\vee\nu_i}}\right)\right)t^{\nu_i}\right)
		E_{\nu_0}\left(c_0\sum_{i=1}^h\left(\kappa_{S_i}^\prime(0)\frac{\theta_i}{t^{\nu_0\vee\nu_i}}
		+o\left(\frac{1}{t^{\nu_0\vee\nu_i}}\right)\right)t^{\nu_0}\right)
	\end{multline*}
	and, by \eqref{eq:ML-asymptotics}, one can check that the limit of the moment generating function (as $t\to\infty$)
	is equal to (we recall that $\kappa_{S_i}^\prime(0)=m_i$)
	\begin{multline*}
		\prod_{i=1}^hE_{\nu_i}\left(c_i\kappa_{S_i}^\prime(0)1_{\{\nu_0\leq\nu_i\}}\theta_i\right)
		E_{\nu_0}\left(c_0\sum_{i=1}^h\kappa_{S_i}^\prime(0)1_{\{\nu_i\leq\nu_0\}}\theta_i\right)\\
		=\prod_{i=1}^h\mathbb{E}\left[e^{c_i\theta_im_i1_{\{\nu_0\leq\nu_i\}}L_i^{\nu_i}(1)}\right]
		\mathbb{E}\left[e^{c_0\sum_{i=1}^h\theta_im_i1_{\{\nu_i\leq\nu_0\}}L_0^{\nu_0}(1)}\right]\\
		=\mathbb{E}\left[e^{\sum_{i=1}^h\theta_i\left(c_im_i1_{\{\nu_0\leq\nu_i\}}L_i^{\nu_i}(1)
		+c_0m_i1_{\{\nu_i\leq\nu_0\}}L_0^{\nu_0}(1)\right)}\right].		
	\end{multline*}
\end{proof}

Now we conclude with the noncentral moderate deviation result.

\begin{proposition}\label{prop:ncMD*}
	Assume that Condition \ref{cond:*} holds.
	Let $\alpha_m(\nu)$ be defined in eq. \eqref{eq:exponents}. Then, for every family of positive numbers 
	$\{a_t:t>0\}$ such that eq. \eqref{eq:MDconditions} holds, the family of random variables $$\left\{\left((ta_t)^{\alpha_m(\nu_0\vee\nu_i)}\frac{S_i(c_iL_i^{\nu_i}(t)+c_0L_0^{\nu_0}(t))}{t}\right)_{i\in\{1,\ldots,h\}}:t>0\right\}$$ 
	satisfies the LDP with speed $1/a_t$ and good rate function $I_{\mathrm{MD}}(\cdot;m)$ defined by
	$$I_{\mathrm{MD}}(x_1,\ldots,x_h;m):=\sup_{(\theta_1,\ldots,\theta_h)\in\mathbb{R}^h}
	\left\{\sum_{i=1}^h\theta_ix_i-\widetilde{\Psi}_m(\theta_1,\ldots,\theta_h)\right\},$$
	where
	\begin{equation}\label{eq:ncMD-function*}
		\widetilde{\Psi}_m(\theta_1,\ldots,\theta_h)=\left\{\begin{array}{ll}
		\sum_{i=1}^h\left(\frac{c_i\kappa_{S_i}^{\prime\prime}(0)}{2}1_{\{\nu_0\leq\nu_i\}}\theta_i^2\right)^{1/\nu_i}\\
		+\left(c_0\sum_{i=1}^h\frac{\kappa_{S_i}^{\prime\prime}(0)}{2}1_{\{\nu_i\leq\nu_0\}}\theta_i^2\right)^{1/\nu_0}&\ \mbox{if}\ m=0\\
		\sum_{i=1}^h\left(c_im_i1_{\{\nu_0\leq\nu_i\}}\theta_i1_{\{\theta_im_i\geq 0\}}\right)^{1/\nu_i}\\
		+\left(c_0\sum_{i=1}^hm_i1_{\{\nu_i\leq\nu_0\}}\theta_i1_{\{\theta_im_i\geq 0\}}\right)^{1/\nu_0}&\ \mbox{if}\ m\neq 0.
		\end{array}\right.
	\end{equation}
\end{proposition}
\begin{proof}
	For every $m\in\mathbb{R}$ we apply the G\"artner Ellis Theorem (Theorem \ref{th:GE}). So we have to show that
	\begin{multline*}
		\lim_{t\to\infty}
		\underbrace{\frac{1}{1/a_t}\log\mathbb{E}\left[e^{\frac{1}{a_t}\sum_{i=1}^h\theta_i
				(ta_t)^{\alpha_m(\nu_0\vee\nu_i)}\frac{S_i(c_iL_i^{\nu_i}(t)+c_0L_0^{\nu_0}(t))}{t}}\right]}_{=:\rho(t)}\\
		=\widetilde{\Psi}_m(\theta_1,\ldots,\theta_h)
		\quad(\mbox{for all}\ (\theta_1,\ldots,\theta_h)\in\mathbb{R}^h),
	\end{multline*}
    where $\widetilde{\Psi}_m(\theta_1,\ldots,\theta_h)$ is defined by eq. \eqref{eq:ncMD-function*}. Indeed, for 
    every $m\in\mathbb{R}^h$, the function $(\theta_1,\ldots,\theta_h)\mapsto\widetilde{\Psi}_m(\theta_1,\ldots,\theta_h)$
    is finite-valued and differentiable (thus the hypotheses of the G\"artner Ellis Theorem are satisfied).
    
    In view of what follows we remark that
	\begin{multline*}
		\rho(t)=a_t\log\mathbb{E}\left[e^{\sum_{i=1}^h\theta_i\frac{S_i(c_iL_i^{\nu_i}(t)+c_0L_0^{\nu_0}(t))}{(ta_t)^{1-\alpha_m(\nu_0\vee\nu_i)}}}
		\right]
		=a_t\log\mathbb{E}\left[e^{\sum_{i=1}^h\kappa_{S_i}(\theta_i/(ta_t)^{1-\alpha_m(\nu_0\vee\nu_i)})(c_iL_i^{\nu_i}(t)+c_0L_0^{\nu_0}(t))}\right]\\
		=a_t\left\{\sum_{i=1}^h\log\mathbb{E}\left[e^{c_i\kappa_{S_i}(\theta_i/(ta_t)^{1-\alpha_m(\nu_0\vee\nu_i)})L_i^{\nu_i}(t)}\right]
		+\log\mathbb{E}\left[e^{c_0\sum_{i=1}^h\kappa_{S_i}(\theta_i/(ta_t)^{1-\alpha_m(\nu_0\vee\nu_i)})L_0^{\nu_0}(t)}\right]\right\}\\
		=a_t\left\{\sum_{i=1}^h\log E_{\nu_i}\left(c_i\kappa_{S_i}(\theta_i/(ta_t)^{1-\alpha_m(\nu_0\vee\nu_i)})t^{\nu_i}\right)
		+\log E_{\nu_0}\left(c_0\sum_{i=1}^h\kappa_{S_i}(\theta_i/(ta_t)^{1-\alpha_m(\nu_0\vee\nu_i)})t^{\nu_0}\right)\right\}.
	\end{multline*}
	
	We start with the case $m=0$. We have
	\begin{multline*}
		\rho(t)=a_t\left\{\sum_{i=1}^h\log E_{\nu_i}\left(\left(\frac{c_i\kappa_{S_i}^{\prime\prime}(0)}{2}\frac{\theta_i^2}{(ta_t)^{\nu_0\vee\nu_i}}
		+o\left(\frac{1}{(ta_t)^{\nu_0\vee\nu_i}}\right)\right)t^{\nu_i}\right)\right.\\
		\left.+\log E_{\nu_0}\left(\sum_{i=1}^h\left(\frac{c_0\kappa_{S_i}^{\prime\prime}(0)}{2}\frac{\theta_i^2}{(ta_t)^{\nu_0\vee\nu_i}}
		+o\left(\frac{1}{(ta_t)^{\nu_0\vee\nu_i}}\right)\right)t^{\nu_0}\right)\right\}
	\end{multline*}
    and, by eq. \eqref{eq:ML-asymptotics}, we obtain
	$$\lim_{t\to\infty}\rho(t)=\sum_{i=1}^h\left(\frac{c_i\kappa_{S_i}^{\prime\prime}(0)}{2}1_{\{\nu_0\leq\nu_i\}}\theta_i^2\right)^{1/\nu_i}
	+\left(c_0\sum_{i=1}^h\frac{\kappa_{S_i}^{\prime\prime}(0)}{2}1_{\{\nu_i\leq\nu_0\}}\theta_i^2\right)^{1/\nu_0}.$$
	Finally the case $m\neq 0$. We have
	\begin{multline*}
        \rho(t)=a_t\left\{\sum_{i=1}^h\log E_{\nu_i}\left(\left(c_i\kappa_{S_i}^\prime(0)\frac{\theta_i}{(ta_t)^{\nu_0\vee\nu_i}}
        +o\left(\frac{1}{(ta_t)^{\nu_0\vee\nu_i}}\right)\right)t^{\nu_i}\right)\right.\\
        \left.+\log E_{\nu_0}\left(\sum_{i=1}^h\left(c_0\kappa_{S_i}^\prime(0)\frac{\theta_i}{(ta_t)^{\nu_0\vee\nu_i}}
        +o\left(\frac{1}{(ta_t)^{\nu_0\vee\nu_i}}\right)\right)t^{\nu_0}\right)\right\},
    \end{multline*}
    and, by eq. \eqref{eq:ML-asymptotics}, we obtain
    $$\lim_{t\to\infty}\rho(t)=\sum_{i=1}^h\left(c_im_i1_{\{\nu_0\leq\nu_i\}}\theta_i1_{\{\theta_im_i\geq 0\}}\right)^{1/\nu_i}\\
    +\left(c_0\sum_{i=1}^hm_i1_{\{\nu_i\leq\nu_0\}}\theta_i1_{\{\theta_im_i\geq 0\}}\right)^{1/\nu_0}.$$
\end{proof}

\begin{remark}\label{rem:*}
	The weak limit in Proposition \ref{prop:weak-convergence*} is a sum of two independent vectors
	$$(B_1,\ldots,B_h)+(\widehat{B}_1,\ldots,\widehat{B}_h)$$
	for some random variables $B_1,\ldots,B_h,\widehat{B}_1,\ldots,\widehat{B}_h$ (say), where 
	$B_1,\ldots,B_h$ are independent, and $\widehat{B}_1,\ldots,\widehat{B}_h$ are conditionally independent given 
	$L_0^{\nu_0}(1)$; actually, if $m\neq 0$, $(\widehat{B}_1,\ldots,\widehat{B}_h)$ is a deterministic vector multiplied 
	by $L_0^{\nu_0}(1)$. This can be related with the expressions of the function $\widetilde{\Psi}_m$ 
	defined by eq. \eqref{eq:ncMD-function*} (for both cases $m=0$ and $m\neq 0$), i.e.
	$$\widetilde{\Psi}_m(\theta_1,\ldots,\theta_h)
	=\sum_{i=1}^hf_i^{1/\nu_i}(\theta_i)+\left(\sum_{i=1}^h\widehat{f}_i(\theta_i)\right)^{1/\nu_0}$$
	for some functions $f_1\ldots,f_h,\widehat{f}_1,\ldots,\widehat{f}_h$ (say).
	
	We also remark that, for every $i\in\{1,\ldots,h\}$, the factors $1_{\{\nu_0\leq\nu_i\}}$ and $1_{\{\nu_i\leq\nu_0\}}$
	appear in $B_i$ and $\widehat{B}_i$ (and also in $f_i$ and $\widehat{f}_i$ as well), respectively. These two factors
	are both equal to 1 if and only if $\nu_0=\nu_i$; otherwise one factor is equal to 1, and the other one is equal to zero.
\end{remark}

\section{Some cases with an explicit expression for $I_{\mathrm{MD}}(\cdot;m)$}\label{sec:explicit-expressions}
In general we do not have an explicit expression for the rate function $I_{\mathrm{LD}}$ concerning the reference LDP (see 
Proposition \ref{prop:LD*}); this is not surprising because we have the same situation for the univariate processes studied
in the literature (see e.g. Proposition 3.1 in \cite{IulianoMacciMeoli}). On the contrary, by taking into account again what
happens for the univariate processes (see e.g. Proposition 3.3 in \cite{IulianoMacciMeoli}), we can wish to get an explicit 
expression for the rate function $I_{\mathrm{MD}}(\cdot;m)$. In this section we present some cases in which this is possible.

In view of what follows it is useful to consider the following notation. 
\begin{itemize}
	\item For $m=0$, we consider the function
	$U_{\nu_i,\kappa_{S_i}^{\prime\prime}(0)}:\mathbb{R}\to[0,\infty]$ defined by
	$$U_{\nu_i,\kappa_{S_i}^{\prime\prime}(0)}(x):=((\nu_i/2)^{\nu_i/(2-\nu_i)}-(\nu_i/2)^{2/(2-\nu_i)})
	\left(\frac{2x^2}{\kappa_{S_i}^{\prime\prime}(0)}\right)^{1/(2-\nu_i)}$$
	for $\kappa_{S_i}^{\prime\prime}(0)>0$, and
	$$U_{\nu_i,\kappa_{S_i}^{\prime\prime}(0)}(x):=\left\{\begin{array}{ll}
		0&\quad\mbox{if}\ x=0\\
		\infty&\quad\mbox{if}\ x\neq 0.
	\end{array}\right.$$
	for $\kappa_{S_i}^{\prime\prime}(0)=0$.
	\item For $m\neq 0$ (note that, in general, we can have either $m_i=0$ or $m_i\neq 0$, but there exists 
	$j\in\{1,\ldots,h\}$ such that $m_j\neq 0$), we consider the function
	$V_{\nu_i,m_i}:\mathbb{R}\to[0,\infty]$ defined by
	$$V_{\nu_i,m_i}(x):=\left\{\begin{array}{ll}
		(\nu_i^{\nu_i/(1-\nu_i)}-\nu_i^{1/(1-\nu_i)})\left(\frac{x_i}{m_i}\right)^{1/(1-\nu_i)}&\quad\mbox{if}\ \frac{x_i}{m_i}\geq 0\\
		\infty&\quad\mbox{if}\ \frac{x_i}{m_i}<0
	\end{array}\right.$$
	for $m_i\neq 0$, and
	$$V_{\nu_i,m_i}(x):=\left\{\begin{array}{ll}
		0&\quad\mbox{if}\ x=0\\
		\infty&\quad\mbox{if}\ x\neq 0
	\end{array}\right.$$
    for $m_i=0$.
\end{itemize}

\subsection{On the case $\nu_0<\min\{\nu_1,\ldots,\nu_h\}$}\label{sub:like-independence}
If $\nu_0<\min\{\nu_1,\ldots,\nu_h\}$, then the function in \eqref{eq:ncMD-function*} reads
$$\widetilde{\Psi}_m(\theta_1,\ldots,\theta_h)=\left\{\begin{array}{ll}
\sum_{i=1}^h\left(\frac{c_i\kappa_{S_i}^{\prime\prime}(0)}{2}\theta_i^2\right)^{1/\nu_i}&\ \mbox{if}\ m=0\\
\sum_{i=1}^h\left(c_im_i\theta_i1_{\{\theta_im_i\geq 0\}}\right)^{1/\nu_i}&\ \mbox{if}\ m\neq 0.
\end{array}\right.$$
So, in both cases $m=0$ and $m\neq 0$, the function $\widetilde{\Psi}_m$ is a sum of $h$ functions, and the $i$-th 
function depends on $\theta_i$ only. In this case $I_{\mathrm{MD}}(\cdot;m)=I_{\mathrm{MD}}(x_1,\ldots,x_h;m)$, 
which is the Legendre-Fenchel transform of $\widetilde{\Psi}_m$, is a sum of $h$ functions, and the $i$-th function 
depends on $x_i$ only; so in some sense we have the \emph{asymptotic independence} of the components. In detail we 
have the following

\begin{proposition}\label{prop:explicit-MD-like-independence*}
	Assume that $\nu_0<\min\{\nu_1,\ldots,\nu_h\}$. Then
	$$I_{\mathrm{MD}}(x_1,\ldots,x_h;m)=\left\{\begin{array}{ll}
		\sum_{i=1}^hU_{\nu_i,c_i\kappa_{S_i}^{\prime\prime}(0)}(x_i)&\ \mbox{if}\ m=0\\
		\sum_{i=1}^hV_{\nu_i,c_im_i}(x_i)&\ \mbox{if}\ m\neq 0.
	\end{array}\right.$$
\end{proposition}
\begin{proof}
	By construction we have
	$$I_{\mathrm{MD}}(x_1,\ldots,x_h;m):=\left\{\begin{array}{ll}
		\sum_{i=1}^h\sup_{\theta_i\in\mathbb{R}}\left\{\theta_ix_i
		-\left(\frac{c_i\kappa_{S_i}^{\prime\prime}(0)}{2}\theta_i^2\right)^{1/\nu_i}\right\}&\ \mbox{if}\ m=0\\
		\sum_{i=1}^h
		\sup_{\theta_i\in\mathbb{R}}\left\{\theta_ix_i
		-\left(c_im_i\theta_i1_{\{\theta_im_i\geq 0\}}\right)^{1/\nu_i}\right\}&\ \mbox{if}\ m\neq 0.
	\end{array}\right.$$

    We start with the case $m=0$. For every $i\in\{1,\ldots,h\}$ we have to check that
    $$\sup_{\theta_i\in\mathbb{R}}\left\{\theta_ix_i
    -\left(\frac{c_i\kappa_{S_i}^{\prime\prime}(0)}{2}\theta_i^2\right)^{1/\nu_i}\right\}
    =U_{\nu_i,c_i\kappa_{S_i}^{\prime\prime}(0)}(x_i).$$
    The case $\kappa_{S_i}^{\prime\prime}(0)=0$ is immediate. For $\kappa_{S_i}^{\prime\prime}(0)>0$ it is easy to 
    check with some computations that the supremum is attained at
    $$\theta_i=\left(\frac{2}{c_i\kappa_{S_i}^{\prime\prime}(0)}\right)^{1/(2-\nu_i)}
    \left(\frac{\nu_i}{2}x_i\right)^{\nu_i/(2-\nu_i)},$$
    and we get
    \begin{multline*}
    	\sup_{\theta_i\in\mathbb{R}}\left\{\theta_ix_i
    	-\left(\frac{c_i\kappa_{S_i}^{\prime\prime}(0)}{2}\theta_i^2\right)^{1/\nu_i}\right\}\\
    	=\left(\frac{2}{c_i\kappa_{S_i}^{\prime\prime}(0)}\right)^{1/(2-\nu_i)}
    	\left(\frac{\nu_i}{2}x_i\right)^{\nu_i/(2-\nu_i)}x_i
    	-\left(\frac{c_i\kappa_{S_i}^{\prime\prime}(0)}{2}\left(\frac{2}{c_i\kappa_{S_i}^{\prime\prime}(0)}\right)^{2/(2-\nu_i)}
    	\left(\frac{\nu_i}{2}x_i\right)^{2\nu_i/(2-\nu_i)}\right)^{1/\nu_i}\\
    	=\left(\frac{2}{c_i\kappa_{S_i}^{\prime\prime}(0)}\right)^{1/(2-\nu_i)}
    	\left(\frac{\nu_i}{2}\right)^{\nu_i/(2-\nu_i)}x_i^{2/(2-\nu_i)}
    	-\left(\frac{2}{c_i\kappa_{S_i}^{\prime\prime}(0)}\right)^{1/(2-\nu_i)}
    	\left(\frac{\nu_i}{2}\right)^{2/(2-\nu_i)}x_i^{2/(2-\nu_i)}\\
    	=\left(\left(\frac{\nu_i}{2}\right)^{\nu_i/(2-\nu_i)}-\left(\frac{\nu_i}{2}\right)^{2/(2-\nu_i)}\right)
    	\left(\frac{2x_i^2}{c_i\kappa_{S_i}^{\prime\prime}(0)}\right)
    	=U_{\nu_i,c_i\kappa_{S_i}^{\prime\prime}(0)}(x_i).
    \end{multline*}
    
    Finally the case $m\neq 0$. For every $i\in\{1,\ldots,h\}$ we have to check that
    $$\sup_{\theta_i\in\mathbb{R}}\left\{\theta_ix_i
    -\left(c_im_i\theta_i1_{\{\theta_im_i\geq 0\}}\right)^{1/\nu_i}\right\}=V_{\nu_i,c_im_i}(x_i).$$
    The case $m_i=0$ is immediate. For $m_i>0$ one can check that: for $x_i<0$ the supremum is equal to infinity
    by taking the limit as $\theta_i\to-\infty$; for $x_i>0$ the supremum is attained at
    $$\theta_i=\frac{(\nu_ix_i)^{\nu_i/(1-\nu_i)}}{(c_im_i)^{1/(1-\nu_i)}}$$
    (after some computations), and we get
    \begin{multline*}
    	\sup_{\theta_i\in\mathbb{R}}\left\{\theta_ix_i
    	-\left(c_im_i\theta_i1_{\{\theta_im_i\geq 0\}}\right)^{1/\nu_i}\right\}\\
    	=\frac{(\nu_ix_i)^{\nu_i/(1-\nu_i)}}{(c_im_i)^{1/(1-\nu_i)}}x_i
    	-\left(c_im_i\frac{(\nu_ix_i)^{\nu_i/(1-\nu_i)}}{(c_im_i)^{1/(1-\nu_i)}}\right)^{1/\nu_i}
    	=\frac{\nu_i^{\nu_i/(1-\nu_i)}}{(c_im_i)^{1/(1-\nu_i)}}x_i^{1/(1-\nu_i)}
    	-\frac{\nu_i^{1/(1-\nu_i)}}{(c_im_i)^{1/(1-\nu_i)}}x_i^{1/(1-\nu_i)}\\
    	=(\nu_i^{\nu_i/(1-\nu_i)}-\nu_i^{1/(1-\nu_i)})\left(\frac{x_i}{c_im_i}\right)^{1/(1-\nu_i)}
    	=V_{\nu_i,c_im_i}(x_i).
    \end{multline*}
    For $m_i<0$ we remark that
    $$\sup_{\theta_i\in\mathbb{R}}\left\{\theta_ix_i
    -\left(c_im_i\theta_i1_{\{\theta_im_i\geq 0\}}\right)^{1/\nu_i}\right\}
    =\sup_{\theta_i\in\mathbb{R}}\left\{-\theta_i(-x_i)
    -\left(c_i(-m_i)(-\theta_i)1_{\{(-\theta_i)(-m_i)\geq 0\}}\right)^{1/\nu_i}\right\}$$
    and, by referring to the previous part of the proof with $-m_i>0$ in place of $m_i$, we get
    $$\sup_{\theta_i\in\mathbb{R}}\left\{-\theta_i(-x_i)
    -\left(c_i(-m_i)(-\theta_i)1_{\{(-\theta_i)(-m_i)\geq 0\}}\right)^{1/\nu_i}\right\}=V_{\nu_i,c_i(-m_i)}(-x_i).$$
    Then, since $V_{\nu_i,c_i(-m_i)}(-x_i)=V_{\nu_i,c_im_i}(x_i)$ by construction, the desired equality holds 
    also for $m_i<0$.
\end{proof}

\begin{remark}\label{rem:asyptotic-independence-only-for-MD*}
	We do not have the asymptotic independence of the components for the rate function 
	$I_{\mathrm{LD}}=I_{\mathrm{LD}}(x_1,\ldots,x_h)$ (as we said above for $I_{\mathrm{MD}}(\cdot;m)$). 
	Indeed, in general (if $\nu_0<\min\{\nu_1,\ldots,\nu_h\}$ holds or not), the function $\Psi$ in Proposition 
	\ref{prop:LD*} cannot be expressed as a sum of $h$ functions such that the $i$-th function depends on $\theta_i$ only.
\end{remark}

\subsection{On the case $\nu_0>\max\{\nu_1,\ldots,\nu_h\}$ with $m\neq 0$}\label{sub:like-a-unique-time-change}
If $\nu_0>\max\{\nu_1,\ldots,\nu_h\}$, then the function in \eqref{eq:ncMD-function*} reads
$$\widetilde{\Psi}_m(\theta_1,\ldots,\theta_h)=\left\{\begin{array}{ll}
\left(c_0\sum_{i=1}^h\frac{\kappa_{S_i}^{\prime\prime}(0)}{2}\theta_i^2\right)^{1/\nu_0}&\ \mbox{if}\ m=0\\
\left(c_0\sum_{i=1}^hm_i\theta_i1_{\{\theta_im_i\geq 0\}}\right)^{1/\nu_0}&\ \mbox{if}\ m\neq 0.
\end{array}\right.$$
Actually here we concentrate our attention on the case $m\neq 0$ only. We have the following

\begin{proposition}\label{prop:prop:explicit-MD-like-a-unique-time-change*}
	Assume that $\nu_0>\max\{\nu_1,\ldots,\nu_h\}$ and $m\neq 0$. Moreover consider the following conditions:\\
	$(\mathrm{i})$ $x_im_i\geq 0$ for every $i\in\{1,\ldots,h\}$;\\
	$(\mathrm{ii})$ if $m_j=0$ some $j\in\{1,\ldots,h\}$, then $x_j=0$.\\
	Then, if we consider the set
	$$\mathcal{A}:=\{(x_1,\ldots,x_h)\in\mathbb{R}^h:(\mathrm{i})\ \mbox{and}\ (\mathrm{ii})\ \mbox{holds}\},$$
	we have
	$$I_{\mathrm{MD}}(x_1,\ldots,x_h;m)=\left\{\begin{array}{ll}
		\max\{V_{\nu_0,c_0m_i}(x_i):i\in\{1,\ldots,h\}\}&\quad (x_1,\ldots,x_h)\in\mathcal{A}\\
		\infty&\quad (x_1,\ldots,x_h)\notin\mathcal{A}.
	\end{array}\right.$$
\end{proposition}
\begin{proof}
	We start noting that, if $h=1$, we can refer to a slight modification of Proposition 3.3 in \cite{IulianoMacciMeoli}
	(here we have to consider $c_0m_1$ in place of $m_1$, and $m_1\neq 0$), and therefore
	$$I_{\mathrm{MD}}(x_1;m_1)=\left\{\begin{array}{ll}
	V_{\nu_0,c_0m_1}(x_1)&\quad x_1m_1\geq 0\\
	\infty&\quad \mbox{otherwise}.
    \end{array}\right.$$	
	So, from now on, we essentially refer to the case $h\geq 2$. We start with $(x_1,\ldots,x_h)\notin\mathcal{A}$, and we
	have two cases.
	\begin{itemize}
		\item If $(\mathrm{i})$ fails, then there exists $j\in\{1,\ldots,h\}$ such that $x_jm_j<0$. So we have two cases:\\
		if $x_j>0$ and $m_j<0$, then we have
		$$\lim_{\theta_j\to+\infty}\theta_jx_j-
		\left(c_0m_j\theta_j1_{\{\theta_jm_j\geq 0\}}\right)^{1/\nu_0}
		=\lim_{\theta_j\to+\infty}\theta_jx_j=+\infty;$$
		if $x_j<0$ and $m_j>0$, then have
		$$\lim_{\theta_j\to-\infty}\theta_jx_j-
		\left(c_0m_j\theta_j1_{\{\theta_jm_j\geq 0\}}\right)^{1/\nu_0}
		=\lim_{\theta_j\to-\infty}\theta_jx_j=+\infty.$$
		\item If $(\mathrm{ii})$ fails, then we have $m_j=0$ some $j\in\{1,\ldots,h\}$ and $x_j\neq 0$; so we immediately get
		$$\sup_{\theta_j\in\mathbb{R}}\left\{\theta_jx_j-
		\left(c_0m_j\theta_j1_{\{\theta_jm_j\geq 0\}}\right)^{1/\nu_0}\right\}
		=\sup_{\theta_j\in\mathbb{R}}\{\theta_jx_j\}=+\infty.$$
	\end{itemize}
    Now we consider $(x_1,\ldots,x_h)\in\mathcal{A}$. By taking into account what we have shown above, we do not lose of
    generality if we assume that $m_i\neq 0$ for every $i\in\{1,\ldots,h\}$; moreover we can say that
    \begin{equation}\label{eq:handled-supremum}
    	 I_{\mathrm{MD}}(x_1,\ldots,x_h;m)=\sup\left\{\sum_{i=1}^h\theta_ix_i
    	-\left(c_0\sum_{i=1}^h\theta_im_i\right)^{1/\nu_0}:\theta_ix_i\geq 0\ \mbox{for every}\ i\in\{1,\ldots,h\}\right\}.
    \end{equation}
    Then we can consider the equations for the stationary points
    \begin{equation}\label{eq:condition-for-stationary-points}
    	x_j=\frac{c_0^{1/\nu_0}}{\nu_0}\left(\sum_{i=1}^h\theta_im_i\right)^{1/\nu_0-1}m_j\quad\mbox{for every}\ j\in\{1,\ldots,h\},
    \end{equation}
    which admit solutions if
    \begin{equation}\label{eq:proportionality-condition}
    	\mbox{there exists $\alpha=\alpha(x_1,\ldots,x_h)\geq 0$ such that $x_i=\alpha m_i$ for every $i\in\{1,\ldots,h\}$.}
    \end{equation}
    We have two cases.
    \begin{enumerate}
    	\item If \eqref{eq:proportionality-condition} holds, then (after some computations, in which we set 
    	$\eta=\sum_{i=1}^h\theta_im_i$) we get
    	\begin{multline*}
    		I_{\mathrm{MD}}(\alpha m_1,\ldots,\alpha m_h;m)=\sup_{\eta\geq 0}\{\alpha\eta-(c_0\eta)^{1/\nu_0}\}\\
    		=\left.\{\alpha\eta-(c_0\eta)^{1/\nu_0}\}\right|_{\eta=(\frac{\alpha\nu_0}{c})^{\nu_0/(1-\nu_0)}}
    		=V_{\nu_0,1}\left(\frac{\alpha}{c_0}\right).
    	\end{multline*}
        Then, in this case, the desired result is proved because we have
        $$V_{\nu_0,c_0m_i}(x_i)=V_{\nu_0,1}\left(\frac{\alpha}{c_0}\right)\quad\mbox{for every}\ i\in\{1,\ldots,h\}.$$
        \item If \eqref{eq:proportionality-condition} fails, then the supremum in \eqref{eq:handled-supremum} is attained when 
        we have $\theta_j=0$ for some $j\in\{1,\ldots,h\}$. Thus
        $$I_{\mathrm{MD}}(x_1,\ldots,x_h;m)=\max\{\varrho_1,\ldots,\varrho_h\},$$
        where
        \begin{multline*}
        	\varrho_j=\sup\left\{\sum_{i=1}^h\theta_ix_i
        	-\left(c_0\sum_{i=1}^h\theta_im_i\right)^{1/\nu_0}:\theta_ix_i\geq 0\ \mbox{for every}\ i\in\{1,\ldots,h\},
        	\ \mbox{and}\ \theta_j=0\right\}
        \end{multline*}
       for every $j\in\{1,\ldots,h\}$. Then, by induction on $h$ (we already know that the case $h=1$ works well), we get
       $$\varrho_j=\max\{V_{\nu_0,c_0m_i}(x_i):i\in\{1,\ldots,h\}\setminus\{j\}\}\quad\mbox{for every}\ j\in\{1,\ldots,h\};$$
       thus we easily obtain
       $$I_{\mathrm{MD}}(x_1,\ldots,x_h;m)=\max\{V_{\nu_0,c_0m_i}(x_i):i\in\{1,\ldots,h\}\}$$
       as desired.
    \end{enumerate}
\end{proof}

\section{Results under Condition \ref{cond:**}}\label{sec:results**}
We start with the reference LDP and the result does not depend on $m$ in Assumption \ref{claim:on-Levy-processes}.
In particular one can check that $J_{\mathrm{LD}}(x_1,\ldots,x_h)=0$ if and only if $(x_1,\ldots,x_h)=(0,\ldots,0)$.

\begin{proposition}\label{prop:LD**}
	Assume that Condition \ref{cond:**} holds.
	Let $\Upsilon$ be the function defined by
	$$\Upsilon(\theta_1,\ldots,\theta_h):=\sum_{j=1}^k\left(\sum_{i=1}^hc_{ij}\kappa_{S_i}(\theta_i)\right)^{1/\nu}
	1_{\{\sum_{i=1}^hc_{ij}\kappa_{S_i}(\theta_i)\geq 0\}}
	\quad(\mbox{for all}\ \theta_1,\ldots,\theta_h\geq 0),$$
	and assume that it is an essentially smooth and lower semicontinuous function. Then
	$$\left\{\left(\frac{S_1\left(\sum_{j=1}^kc_{1j}L_j^\nu(t)\right)}{t},\ldots,
	\frac{S_h\left(\sum_{j=1}^kc_{hj}L_j^\nu(t)\right)}{t}\right):t>0\right\}$$
	satisfies the LDP with speed function $t$ and rate function $J_{\mathrm{LD}}$ defined by
	$$J_{\mathrm{LD}}(x_1,\ldots,x_h):=\sup_{(\theta_1,\ldots,\theta_h)\in\mathbb{R}^h}\left\{\sum_{i=1}^h\theta_ix_i
	-\Upsilon(\theta_1,\ldots,\theta_h)\right\}\quad(\mbox{for all}\ (x_1,\ldots,x_h)\in\mathbb{R}^h.$$
\end{proposition}
\begin{proof}
	The desired LDP can be derived by applying the G\"artner Ellis Theorem (i.e. Theorem \ref{th:GE}). In fact, for all
	$t>0$ (and for all $(\theta_1,\ldots,\theta_h)\in\mathbb{R}^h$), we have
	\begin{multline*}
		\mathbb{E}\left[e^{\sum_{i=1}^h\theta_iS_i\left(\sum_{j=1}^kc_{ij}L_j^\nu(t)\right)}\right]
		=\mathbb{E}\left[\mathbb{E}\left[e^{\sum_{i=1}^h\theta_iS_i\left(\sum_{j=1}^kc_{ij}r_j\right)}\right]_{r_1=L_1^\nu(t),
			\ldots,r_k=L_k^\nu(t)}\right]\\
		=\mathbb{E}\left[\prod_{i=1}^h\mathbb{E}\left[e^{\theta_iS_i(1)}\right]^{\sum_{j=1}^kc_{ij}L_j^\nu(t)}\right]
		=\mathbb{E}\left[e^{\sum_{i=1}^h(\sum_{j=1}^kc_{ij}L_j^\nu(t))\kappa_{S_i}(\theta_i)}\right]\\
		=\mathbb{E}\left[e^{\sum_{j=1}^k(\sum_{i=1}^hc_{ij}\kappa_{S_i}(\theta_i))L_j^\nu(t)}\right]
	\end{multline*}
	whence we obtain (by the independence of the involved inverse stable subordinators, and by the expressions of the
	involved moment generating functions given by \eqref{eq:MGF-inverse-stable-sub})
	\begin{multline*}
		\frac{1}{t}\log\mathbb{E}\left[e^{\sum_{i=1}^h\theta_iS_i\left(\sum_{j=1}^kc_{ij}L_j^\nu(t)\right)}\right]\\
		=\sum_{j=1}^k\frac{1}{t}\log\mathbb{E}\left[e^{\sum_{i=1}^hc_{ij}\kappa_{S_i}(\theta_i)L_j^\nu(t)}\right]
		=\sum_{j=1}^k\frac{1}{t}\log E_\nu\left(\sum_{i=1}^hc_{ij}\kappa_{S_i}(\theta_i)t^\nu\right).
	\end{multline*}
	Finally, if we take the limit as $t$ tends to infinity, we get
	$$\lim_{t\to\infty}\frac{1}{t}\log\mathbb{E}\left[e^{\sum_{i=1}^h\theta_iS_i\left(\sum_{j=1}^kc_{ij}L_j^\nu(t)\right)}\right]
	=\Upsilon(\theta_1,\ldots,\theta_h)$$
	by taking into account eq. \eqref{eq:ML-asymptotics} with $\alpha=\nu$.
\end{proof}

As happens for the analogue results in Section \ref{sec:results*}, in the following results we have to distinguish the cases 
$m=0$ and $m\neq 0$, where $0\in\mathbb{R}^h$ is the null vector. We start with the weak convergence result and, again, we take 
the limit of moment generating functions (instead of characteristic functions). 

\begin{proposition}\label{prop:weak-convergence**}
	Assume that Condition \ref{cond:**} holds.	
	Let $\alpha_m(\nu)$ be defined in eq. \eqref{eq:exponents}. We have the following statements.
	\begin{itemize}
		\item If $m=0$, then $\left\{\left(t^{\alpha_m(\nu)}
		\frac{S_i\left(\sum_{j=1}^kc_{ij}L_j^\nu(t)\right)}{t}\right)_{i\in\{1,\ldots,h\}}:t>0\right\}$ 
		converges weakly to
		$$\left(\sqrt{\sum_{j=1}^kc_{ij}L_j^\nu(1)\kappa_{S_i}^{\prime\prime}(0)}Z_i\right)_{i\in\{1,\ldots,h\}},$$
		where $Z_1,\ldots,Z_h$ are independent standard Normal distributed random variables, and independent of all the other 
		random variables.
		\item If $m\neq 0$, then $\left\{\left(t^{\alpha_m(\nu)}
		\frac{S_i\left(\sum_{j=1}^kc_{ij}L_j^\nu(t)\right)}{t}\right)_{i\in\{1,\ldots,h\}}:t>0\right\}$ 
		converges weakly to
		$$\left(m_i\sum_{j=1}^kc_{ij}L_j^\nu(1)\right)_{i\in\{1,\ldots,h\}}.$$
	\end{itemize}
\end{proposition}
\begin{proof}
	In both cases $m=0$ and $m\neq 0$ we study suitable limits (as $t\to\infty$) of suitable moment generating functions,
	which can be expressed in terms of Mittag-Leffler functions (see \eqref{eq:MGF-inverse-stable-sub}). Moreover
	\begin{multline*}
		\mathbb{E}\left[e^{\sum_{i=1}^h\theta_it^{\alpha_m(\nu)}\frac{S_i\left(\sum_{j=1}^kc_{ij}L_j^\nu(t)\right)}{t}}\right]
		=\mathbb{E}\left[e^{\sum_{i=1}^h\theta_i\frac{S_i\left(\sum_{j=1}^kc_{ij}L_j^\nu(t)\right)}{t^{1-\alpha_m(\nu)}}}\right]\\
		=\mathbb{E}\left[\mathbb{E}\left[e^{\sum_{i=1}^h\theta_i\frac{S_i\left(\sum_{j=1}^kc_{ij}r_j\right)}{t^{1-\alpha_m(\nu)}}}
		\right]_{r_1=L_1^\nu(t),\ldots,r_k=L_k^\nu(t)}\right]
		=\mathbb{E}\left[e^{\sum_{i=1}^h\left(\sum_{j=1}^kc_{ij}L_j^\nu(t)\right)\kappa_{S_i}(\theta_i/t^{1-\alpha_m(\nu)})}\right]\\
		=\prod_{j=1}^k\mathbb{E}\left[e^{\sum_{i=1}^hc_{ij}\kappa_{S_i}(\theta_i/t^{1-\alpha_m(\nu)})L_j^\nu(t)}\right]
		=\prod_{j=1}^kE_\nu\left(\sum_{i=1}^hc_{ij}\kappa_{S_i}(\theta_i/t^{1-\alpha_m(\nu)})t^\nu\right).
	\end{multline*}

	We start with the case $m=0$. We have
	$$\mathbb{E}\left[e^{\sum_{i=1}^h\theta_it^{\alpha_m(\nu)}\frac{S_i\left(\sum_{j=1}^kc_{ij}L_j^\nu(t)\right)}{t}}\right]
	=\prod_{j=1}^kE_\nu\left(\sum_{i=1}^hc_{ij}\left(\frac{\kappa_{S_i}^{\prime\prime}(0)}{2}\frac{\theta_i^2}{t^\nu}
	+o\left(\frac{1}{t^\nu}\right)\right)t^\nu\right)$$
	and, by eq. \eqref{eq:ML-asymptotics}, we get
	$$\lim_{t\to\infty}\mathbb{E}\left[e^{\sum_{i=1}^h\theta_it^{\alpha_m(\nu)}\frac{S_i\left(\sum_{j=1}^kc_{ij}L_j^\nu(t)\right)}{t}}\right]
	=\prod_{j=1}^kE_\nu\left(\sum_{i=1}^h\frac{c_{ij}\kappa_{S_i}^{\prime\prime}(0)}{2}\theta_i^2\right)
	=\mathbb{E}\left[e^{\sum_{i=1}^h\theta_i\sqrt{\sum_{j=1}^kc_{ij}L_j^\nu(1)\kappa_{S_i}^{\prime\prime}(0)}Z_i}\right]$$
	(the second equality can be easily checked).
	
	Finally the case $m\neq 0$. We have 
	$$\mathbb{E}\left[e^{\sum_{i=1}^h\theta_it^{\alpha_m(\nu)}\frac{S_i\left(\sum_{j=1}^kc_{ij}L_j^\nu(t)\right)}{t}}\right]
	=\prod_{j=1}^kE_\nu\left(\sum_{i=1}^hc_{ij}\left(\kappa_{S_i}^\prime(0)\frac{\theta_i}{t^\nu}
	+o\left(\frac{1}{t^\nu}\right)\right)t^\nu\right)$$
	and, by eq. \eqref{eq:ML-asymptotics}, we get
	$$\lim_{t\to\infty}\mathbb{E}\left[e^{\sum_{i=1}^h\theta_it^{\alpha_m(\nu)}\frac{S_i\left(\sum_{j=1}^kc_{ij}L_j^\nu(t)\right)}{t}}\right]
	=\prod_{j=1}^kE_\nu\left(\sum_{i=1}^hc_{ij}\kappa_{S_i}^\prime(0)\theta_i\right)
	=\mathbb{E}\left[e^{\sum_{i=1}^h\theta_im_i\sum_{j=1}^kc_{ij}L_j^\nu(1)}\right]$$
	(the second equality can be easily checked; we recall that $\kappa_{S_i}^\prime(0)=m_i$).
\end{proof}

Now we conclude with the noncentral moderate deviation result.

\begin{proposition}\label{prop:ncMD**}
	Assume that Condition \ref{cond:**} holds.
	Let $\alpha_m(\nu)$ be defined in eq. \eqref{eq:exponents}. Then, for every family of positive numbers 
	$\{a_t:t>0\}$ such that eq. \eqref{eq:MDconditions} holds, the family of random variables
	$$\left\{\left((ta_t)^{\alpha_m(\nu)}
	\frac{S_i\left(\sum_{j=1}^kc_{ij}L_j^\nu(t)\right)}{t}\right)_{i\in\{1,\ldots,h\}}:t>0\right\}$$ 
	satisfies the LDP with speed $1/a_t$ and good rate function $J_{\mathrm{MD}}(\cdot;m)$ defined by
	$$J_{\mathrm{MD}}(x_1,\ldots,x_h;m):=\sup_{(\theta_1,\ldots,\theta_h)\in\mathbb{R}^h}
	\left\{\sum_{i=1}^h\theta_ix_i-\widetilde{\Upsilon}_m(\theta_1,\ldots,\theta_h)\right\},$$
	where
	\begin{equation}\label{eq:ncMD-function**}
		\widetilde{\Upsilon}_m(\theta_1,\ldots,\theta_h)=\left\{\begin{array}{ll}
			\sum_{j=1}^k\left(\sum_{i=1}^h\frac{c_{ij}\kappa_{S_i}^{\prime\prime}(0)}{2}\theta_i^2\right)^{1/\nu}&\ \mbox{if}\ m=0\\
			\sum_{j=1}^k\left(\sum_{i=1}^hc_{ij}m_i\theta_i1_{\{\theta_im_i\geq 0\}}\right)^{1/\nu}&\ \mbox{if}\ m\neq 0.
		\end{array}\right.
	\end{equation}
\end{proposition}
\begin{proof}
	As in the proof of Proposition \ref{prop:ncMD*}, 
	for every $m\in\mathbb{R}$ we apply the G\"artner Ellis Theorem (Theorem \ref{th:GE}). So we have to show that
	\begin{multline*}
		\lim_{t\to\infty}
		\underbrace{\frac{1}{1/a_t}\log\mathbb{E}\left[e^{\frac{1}{a_t}\sum_{i=1}^h\theta_i(ta_t)^{\alpha_m(\nu)}
			\frac{S_i\left(\sum_{j=1}^kc_{ij}L_j^\nu(t)\right)}{t}}\right]}_{=:\rho(t)}\\
		=\widetilde{\Upsilon}_m(\theta_1,\ldots,\theta_h)
		\quad(\mbox{for all}\ (\theta_1,\ldots,\theta_h)\in\mathbb{R}^h),
	\end{multline*}
	where $\widetilde{\Upsilon}_m(\theta_1,\ldots,\theta_h)$ is defined by eq. \eqref{eq:ncMD-function**}. Indeed, for 
	every $m\in\mathbb{R}^h$, the function $(\theta_1,\ldots,\theta_h)\mapsto\widetilde{\Upsilon}_m(\theta_1,\ldots,\theta_h)$
	is finite-valued and differentiable (thus the hypotheses of the G\"artner Ellis Theorem are satisfied).
	
	In view of what follows we remark that
	\begin{multline*}
		\rho(t)=a_t\log\mathbb{E}\left[e^{\sum_{i=1}^h\theta_i\frac{S_i\left(\sum_{j=1}^kc_{ij}L_j^\nu(t)\right)}
			{(ta_t)^{1-\alpha_m(\nu)}}}\right]
		=a_t\log\mathbb{E}\left[e^{\sum_{i=1}^h\left(\sum_{j=1}^kc_{ij}L_j^\nu(t)\right)\kappa_{S_i}(\theta_i/(ta_t)^{1-\alpha_m(\nu)})}\right]\\
		=a_t\log\mathbb{E}\left[e^{\sum_{j=1}^k\left(\sum_{i=1}^hc_{ij}\kappa_{S_i}(\theta_i/(ta_t)^{1-\alpha_m(\nu)})\right)L_j^\nu(t)}\right]
		=a_t\sum_{j=1}^k\log\mathbb{E}\left[e^{\sum_{i=1}^hc_{ij}\kappa_{S_i}(\theta_i/(ta_t)^{1-\alpha_m(\nu)})L_j^\nu(t)}\right]\\
		=a_t\sum_{j=1}^k\log E_\nu\left(\sum_{i=1}^hc_{ij}\kappa_{S_i}(\theta_i/(ta_t)^{1-\alpha_m(\nu)})t^\nu\right).
	\end{multline*}
	
	We start with the case $m=0$. We have
	$$\rho(t)=a_t\sum_{j=1}^k\log E_\nu\left(\sum_{i=1}^hc_{ij}\left(\frac{\kappa_{S_i}^{\prime\prime}(0)}{2}\frac{\theta_i^2}{(ta_t)^\nu}
	+o\left(\frac{1}{(ta_t)^\nu}\right)\right)t^\nu\right)$$
	and, by eq. \eqref{eq:ML-asymptotics}, we obtain
	$$\lim_{t\to\infty}\rho(t)=\sum_{j=1}^k\left(\sum_{i=1}^hc_{ij}\frac{\kappa_{S_i}^{\prime\prime}(0)}{2}\theta_i^2\right)^{1/\nu}.$$
	Finally the case $m\neq 0$. We have
	$$\rho(t)=a_t\sum_{j=1}^k\log E_\nu\left(\sum_{i=1}^hc_{ij}\left(\kappa_{S_i}^\prime(0)\frac{\theta_i}{(ta_t)^\nu}
	+o\left(\frac{1}{(ta_t)^\nu}\right)\right)t^\nu\right)$$
	and, by eq. \eqref{eq:ML-asymptotics}, we obtain (we recall that $\kappa_{S_i}^\prime(0)=m_i$)
	$$\lim_{t\to\infty}\rho(t)=\sum_{j=1}^k\left(\sum_{i=1}^hc_{ij}m_i\theta_i1_{\{\theta_im_i\geq 0\}}\right)^{1/\nu}.$$
\end{proof}

\paragraph{Acknowledgements.}
The authors thank Carlo Sinestrari for the final part of the proof of Proposition 
\ref{prop:prop:explicit-MD-like-a-unique-time-change*} (in which it is assumed that 
\eqref{eq:proportionality-condition} fails), and two referees for their useful comments.


\begin{thebibliography}{sp}
\bibitem{BeghinMacciSPL2022}
L. Beghin, C. Macci (2022) Non-central moderate deviations for compound fractional
Poisson processes. Statist. Probab. Lett. 185, Paper No. 109424, 8 pp. (2022).
\bibitem{CD}
F. Colantoni, M. D'Ovidio (2023) On the inverse gamma subordinator. Stochastic Anal. 
Appl. 41, 999--1024.
\bibitem{DemboZeitouni}
A. Dembo, O. Zeitouni (1998) Large Deviations Techniques and Applications (Second 
Edition). Springer-Verlag, New York.
\bibitem{GiulianoMacci}
R. Giuliano, C. Macci (2023) Some examples of noncentral moderate deviations for 
sequences of real random variables. Mod. Stoch. Theory Appl. 10, 111--144.
\bibitem{GorenfloKilbasMainardiRogosin}
R. Gorenflo, A.A. Kilbas, F. Mainardi, S.V. Rogosin (2014) 
Mittag-Leffler Functions, Related Topics and Applications. Springer,
New York.
\bibitem{GKL}
N. Gupta, A. Kumar, N. Leonenko (2021) Stochastic models with mixtures of tempered
stable subordinators. Math. Commun. 26, 77--99.
\bibitem{IulianoMacciMeoli}
A. Iuliano, C. Macci, A. Meoli (2025) Noncentral moderate deviations for 
time-changed Lévy processes with inverse of stable subordinators.
Mod. Stoch. Theory Appl. 12, 203--224.
\bibitem{KV}
A. Kumar, P. Vellaisamy (2015) Inverse tempered stable subordinators. Stat. Probab.
Lett. 103, 134--141.
\bibitem{LeeMacci}
J. Lee, C. Macci (2024) Noncentral moderate deviations for fractional
Skellam processes. Mod. Stoch. Theory Appl. 11, 43--61.
\bibitem{LeonenkoMacciPacchiarotti}
N. Leonenko, C. Macci, B. Pacchiarotti (2021), Large deviations for a class of 
tempered subordinators and their inverse processes. Proc. Roy. Soc. Edinburgh 
Sect. A 151, 2030--2050.
\bibitem{LMSS}
N.N. Leonenko, M.M. Meerschaert, R.L. Schilling, A. Sikorskii (2014), Correlation 
structure of time-changed Lévy processes.
Commun. Appl. Ind. Math. 6, No. 1, Article ID 483, 22 p. 
\bibitem{MS}
M.M. Meerschaert, P. Straka (2013) Inverse stable subordinators.
Math. Model. Nat. Phenom. 8, 1--16.

\end{thebibliography}
\end{document}